\input amstex
\documentstyle{amsppt}

\magnification=\magstep1
\hsize=14 truecm
\vsize=19.1 truecm
\hoffset=+0.5 truecm
\voffset=+2 truecm
\headline{\hfill}

\centerline {\bf TOWARDS A HALPHEN THEORY}
\centerline {\bf OF LINEAR SERIES ON CURVES}
\medskip
\centerline {by}
\medskip
\centerline { L. CHIANTINI}
\centerline {Universit\`a di Siena}
\centerline {Dipartimento di Matematica}
\centerline {Via del Capitano 15}
\centerline {53100 Siena (Italy)}
\medskip
\centerline {and}
\medskip
\centerline { C. CILIBERTO} 
\centerline {Universit\`a di Roma "Tor Vergata"}
\centerline {Dipartimento di Matematica}
\centerline {Via dellla Ricerca Scientifica}
\centerline {00133 Roma (Italy)}
\bigskip

\topmatter

\abstract 
A linear series $g^N_\delta$ on a curve $C\subset\bold P^3$ is {\it primary} when it does not contain the series cut by
planes. For such series, we provide a lower bound for the degree $\delta$, in terms of deg($C$), g($C$) and of the number $s=\min\{i:h^0\Cal I_C(i)\neq 0\}$. Examples show that the bound is sharp.\par
Extensions to the case of general linear series and to the case of curves in higher projective spaces are considered.
\endabstract

\endtopmatter

\document
\bigskip

\heading 0. Introduction. \endheading

\smallskip

Let $C$ be a smooth, complete curve of genus $g$ over an algebraically closed field 
$k$ of characteristic zero. If $C$ has a line bundle ${\Cal L}$ generated by
global sections, of degree $d$ with $h^0{\Cal L}=r+1$, i.e. a complete, base point free
$g^r_d$, determining a morphism of $C$ to ${\bold P}^r$ which is birational onto
its image $\Gamma$, then the classical theory of Castelnuovo provides a bound for the
genus $g$ of $C$ in terms of $d$ and $r$. Indeed one has\par

$$g\leq p_a(\Gamma)\leq G(r;d):= {\binom m 2}(r-1)+m\epsilon\eqno (!)$$

\noindent where $d=m(r-1)+\epsilon$ with $0\leq \epsilon\leq r-2$ and this bound is sharp
(see [H2]). \par

A natural question in this circle of ideas is to look for a similar sharp upper bound
$G(r_1,...,r_n;d_1,...,d_n)$ for the genus of a curve $C$ having several different
linear series $g^{r_i}_{d_i}$, $i=1,...,n$. This question has been considered by Accola
in [A] and it can be solved, in some cases, with the same techniques used to prove
Castelnuovo's bound (!).\par

Another classical direction in which Castelnuovo's theory has been developed, has been the
search for a sharp upper bound $G(r;d,s)$ for the genus of a non degenerate curve of
degree $d$ in ${\bold P}^r$ not lying on any irreducible surface $S$ of degree $\sigma<
s$. This is the so-called {\it Halphen problem} inasmuch as it was first considered by
Halphen in a famous memoir of 1882 (see [Ha], and, for more recent references, [H2] and
[CCD1]). Therefore, taking into account Accola's viewpoint mentioned above, a natural
development of both Castelnuovo's and Halphen's theories consists in looking for similar
bounds for the genus of non degenerate curves of degree $d$ in ${\bold P}^r$ not lying on
any irreducible surface $S$ of degree $\sigma< s$ and having the further property of
possessing an additional linear series $g^N_\delta$, besides the $g^r_d$ cut out by the
hyperplanes of ${\bold P}^r$. This is the question we deal with in the present paper,
which is essentially devoted to set the foundations of this Castelnuovo-Halphen theory
for curves with {\it extra series} and to study a few relevant examples. Hence this note
can be seen as a continuation of the mentioned papers [A], [CCD1] and [H2] as well as of
other papers like [B], [C1], [CK], [CL], [L], [M], [P], ecc. in which questions like the
ones we consider here are dealt with in some particular instances.\par
We recall that the classical Castelnuovo's (Halphen's) theory is based on the idea of
bounding from below the Hilbert function $h_\Gamma(n)$ of an irreducible, non degenerate
curve $\Gamma$ in ${\bold P}^r$ (not lying on any surface of degree $\sigma<s$), or, if
one wishes, in bounding from below the function $h^0{\Cal O}_C(nH)$, where $C$ is, as
above, the normalization of $\Gamma$ and $H$ is the pull-back to $C$ of a hyperplane in 
${\bold P}^r$. In the same vein, our approach here consists in trying to bound from
below the {\it twisted Hilbert function} $h^0{\Cal O}_C(D+nH)$, where $|D|$
 is the extra linear series on $C$. Our main results in this direction are contained in
\S 2, while \S 1 is devoted to generalities and to preparatory material. In particular
we note our {\it descent lemmas} 2.3 and 2.4 which play an essential role in our
approach and whose proofs are consequences of previous results we proved
(in collaboration with Di Gennaro) in [CCD1] and [CCD2]. As in the classical theory,
this sort of descent lemma would be almost useless if not accompanied by a numerical
analysis aimed to find the {\it optimal function} minimizing the twisted Hilbert
functions of the curves we deal with. This numerical analysis is carried out in \S 3
and, as sometimes happens also in the classical case (see [CCD2]), it is the hardest part
of the whole story. Next sections are devoted to show some applications. First we concentrate 
in \S 4 on space curves, i.e. curves in ${\bold P}^3$, to which theorem 2.3 applies. 
Our main results (see proposition 4.1 and theorem 4.4) give  a lower bound for the degree 
$\delta$ of a $g^N_\delta$ on a curve of degree $d$ and genus $g$ not lying on any surface of 
degree $\sigma<s$ in ${\bold P}^3$, under the additional hypothesis that the series is {\it 
primary}, i.e. it does not contain the series cut out by the planes. The bound is in terms of 
$d, g$ and $s$ and it is sharp. In particular we are able to bound the gonality and 
the degree of a linear series of dimension $2$, on curves of certain types 
like complete intersections, Halphen curves, curves on quadrics etc. and to 
classify those curves for which the bound is attained. As we indicated before, 
some results of this type can already be found in the papers we mentioned above. 
We are also able to give information for the non-primary
series, but in this case the numerical analysis needed for producing sharp 
bounds, becomes  more complicated, and will not be developed here. 
Finally in \S 5, using proposition 2.4, we make some application to curves in higher projective
spaces. Our results in this case are not sharp as for space curves, however we are
able to provide another proof of Accola's bound for the gonality of Castelnuovo's curves
and we give also a lower bound for the degree of a $g^N_\delta$ with $N\geq 2$ on such
curves; this bound is sharp in some numerical range for $d,N$.\par

We believe that the method presented here for the study of extra linear series can be usefully applied in many other circumstances: it will be enough to use
the results of \S 2 with an appropriate numerical analysis.  
This paper is essentially devoted to introduce the method, show its first application and
prove a few results for some specific classes of curves most of which are new, and, we feel, not
without some interest.\par

\heading 1. Notation. \endheading
\smallskip

Let $\Gamma\subset\bold P^r$ be a projective irreducible non degenerate curve of degree $d$ 
and let $\sigma: C\to\Gamma$ be its normalization. Let $g$ be the genus of $C$. Let $\omega_C$
be the canonical sheaf on  $C$ and let $H$ be the pull-back on $C$ of the hyperplane divisor of
$\Gamma$.\par 

For any divisor $D$ on $C$, we denote by $\Cal O_C(D)$ the sheaf of sections of the
line bundle associated to $D$ and $\delta$ will denote the degree of $D$.\par\medskip

\proclaim{Definition 1.1} The {\rm Hilbert function of $D$ with respect to $H$} is the
function
$$\phi_D(n) = h^0\Cal O_C(D+nH). $$ \par
We say that $D$ is {\rm $H$-primary} if it is effective and $h^0\Cal O_C(D-H)=0$, i.e.
if $\phi_D(n)=0$ for $n<0$ while  $\phi_D(n)>0$ when $n\ge 0$.  \par
$\Delta_D^1$ and $\Delta_D^2$ will denote the first and the second diffence of $\phi_D$, e.g.
$$\gather \Delta^1_D(n)=\phi_D(n)-\phi_D(n-1) \\
 \Delta^2_D(n) = \Delta^1_D(n)-\Delta^1_D(n-1) = \phi_D(n)-2\phi_D(n-1)+\phi_D(n-2). 
\endgather$$
Clearly  $$\phi_D(n) = \sum_{i=-\infty}^n\sum_{j=-\infty}^i\Delta^2_D(j) $$
and for $n\ll 0$ we have $\phi_D(n)=\Delta^1_D(n)=\Delta^2_D(n)=0$.
\endproclaim

It follows by the definition that $\phi_D(n+1)>\phi_D(n)$ whenever $\phi_D(n+1)>0$.
\par
When $\Gamma$ is smooth and arithmetically normal, then taking $D=0$, $\phi_D$ turns
out to be the usual Hilbert function of $\Gamma=C$.

\proclaim{Definition 1.2} Define the {\rm level (of speciality) of $D$ with respect to $H$}
as the number $$e(D)=\max\{n:h^1\Cal O_C(D+nH)>0\}.$$
The number $e:=e(0)$ is also called the {\rm level of $C$}. The divisor $D$ is {\rm non-special} if and only if $e(D)<0$. If $D$ is effective, then $e(D)\le e$.
\endproclaim

\proclaim{Remark 1.3} \rm If $D$ is $H-$primary, then 

$$\phi_D(0)=\Delta^1_D(0)
=\Delta^2_D(0)= h^0\Cal O_C(D).$$\par\smallskip
For $n>e(D)$, by Riemann-Roch we have $\phi_D(n)=nd+\delta-g+1$.\par\smallskip
For $n>e(D)+1$ we have $\Delta^1_D(n)=d$.\par\smallskip
For $n>e(D)+2$ we have $\Delta^2_D(n)=0$.\par\smallskip
\endproclaim

\noindent Conversely, we have:

\proclaim{Proposition 1.4} Assume $\phi_D(n)=nd+\delta-g+1$. Then $n>e(D)$.\par
Assume $\Delta^1_D(m)=d$ for all $m\ge n$. Then $n>e(D)+1$.\par \smallskip
Assume $\Delta^2_D(m)=0$ for all $m\ge n$. Then $n>e(D)+2$.
\endproclaim
\demo{Proof} The first assertion follows by Riemann-Roch. The second and the third
are easy consequences of the first.
$\diamondsuit$\enddemo

\proclaim{Proposition 1.5}  Put $\Cal O_C(D')\simeq\omega_C\otimes\Cal O_C(-D)$. Then
$$ \Delta^2_D(i)=\Delta^2_{D'}(2-i) $$
for all $i$.\endproclaim
\demo{Proof} Follows  by Serre duality.$\diamondsuit$\enddemo

\proclaim{Proposition 1.6}  $\Delta^1_D$ is non decreasing, i.e. $\Delta^2_D(n)\ge 0$
for all $n$.\endproclaim
\demo{Proof} Let $Z$ be a general hyperplane section of $C$ and look at the exact  commutative diagram:
$$\CD  0 @. 0 @. \\@VVV  @VVV @. \\
0 \to H^0\Cal O_C(D+(n-2)H) @>>> H^0\Cal O_C(D+(n-1)H) @>\alpha_{n-1}>> H^0\Cal O_Z \\
        @VVV     @VVV   @| \\
0 \to H^0\Cal O_C(D+(n-1)H) @>>> H^0\Cal O_C(D+nH) @>\alpha_{n}>> H^0\Cal O_Z \\
\endCD   $$
in which $\Delta^1_D(n-1)=$rank$(\alpha_{n-1})$ and $\Delta^1_D(n)=$rank$(\alpha_n)$.
$\diamondsuit$\enddemo

\bigskip
\heading  2. Basic Results. \endheading
\smallskip

We collect here some basic lemmas  which will be necessary for our numerical analysis
of the Hilbert function of $D$.\par
While the speciality lemma is well-known, the descent lemmas are  new and useful.

\proclaim{Proposition 2.1 (Speciality Lemma)} Let $\Gamma$ be a curve in $\bold P^r$, $\sigma:C\to
\Gamma$ its normalization, of genus $g$, let $|D|$ be a base point free
linear series on $C$ and put $N+1=h^0\Cal O_C(D)$. Then  $e(D)\le e \le \delta-N-1$.\par
If $e=\delta-N-1$, then either $e\leq 0$ and $g\leq 1$, or the image on $\Gamma$ of the general 
divisor of $|D|$ lies on a line and either $N=1$ or $N=r=2$.\par
If $e=\delta-N-2$, then either $e\leq 0$, or the image on $\Gamma$ of the general 
divisor of $|D|$ lies on a plane.\par
Furthermore, if $D_0\in |D|$ is a general  divisor and $\Cal I$ is its ideal sheaf in the
 projective space $\bold P^r$, then $h^1\Cal I(e)\ge N$. 
\endproclaim
\demo{Proof} Since $|D|$ is base point free, we may assume $D_0\in  |D|$ disjoint
from the pull back of the singular locus of $\Gamma$. \par
The exact sequence $0\to \Cal I_{D_0,C}(\omega_C)\to\Cal O_C(\omega_C)\to
\Cal O_{D_0}\to 0$ gives in cohomology:
$$ H^0\Cal O_C(\omega_C)\to H^0\Cal O_{D_0}\to H^1\Cal I_{D_0,C}(\omega_C)\to H^1\Cal 
O_C(\omega_C)\to 0$$
where $h^1\Cal I_{D_0,C}(\omega_C)=h^1\Cal O_C(\omega_C(-D))=h^0\Cal O_C(D)=N+1$, so that the image of the first map has dimension $\delta-N$.\par
On the other hand, since $h^0\omega_C(-eH)>0$, we have a factorization:
$$ \CD  H^0\Cal O_C(\omega_C) @>>>  H^0\Cal O_{D_0} \\
           @AAA  @| @. \\
        H^0\Cal O_{\bold P^r}(e)   @>\beta>>  H^0\Cal O_{D_0}    \endCD $$
\noindent so that dim(Im$(\beta))\le \delta-N$, which means $h^1\Cal I(e)\ge N$.\par
Finally, for any set of distinct points in $\bold P^r$, we  have  
dim(Im($\beta))\ge e+1$, and, if $e\geq 1$, equality holds if and only if $D_0$ lies on a line,
while dim(Im($\beta))=e+2$ implies that $D_0$ is contained on a plane. 
$\diamondsuit$\enddemo

More generally, Im$(\beta)$ has dimension $\gg e$, i.e. $\delta-N-1\gg e$, unless the points of $D_0$ belong to some linear subspace.\par
 In particular:

\proclaim{Corollary 2.2} If $C$ is hyperelliptic then $e=0$.\par
If $C$ is trigonal, then $e\le 1$ and $e=1$ if and only if  the image on $\Gamma$ of the general 
divisor of the $g^1_3$ spans a line.\par
If $C$ is quadrigonal, then $e\le 2$ and $e=2$ implies that the image on $\Gamma$ of the general 
divisor of the  $g^1_4$ spans a line, whereas $e=1 $ implies that the image on $\Gamma$ of the 
general  divisor of the  $g^1_4$ spans a line or a plane. 
\endproclaim

 Our next result only works in $\bold P^3$. It provides our main tool for 
handling the function $\Delta^2_D$.

\proclaim{Theorem 2.3 (Descent Lemma)} Let $\Gamma$ be a curve in $\bold P^3$, 
$\sigma:C\to \Gamma$ its normalization, $D$ a divisor on $C$. \par
If $t$ is any number which satisfies: 
$$ \binom{t-1}2< \Delta^1_D(n)<d-\binom t2$$
\noindent and $h^0\Cal I_\Gamma(t-1)=0$, then we have $\Delta^2_D(n+1)\ge t$.
\endproclaim
\demo{Proof} For a general point $\eta\in(\bold P^3)^*$ call $\pi_\eta$ the corresponding
plane and let $Z_\eta$ be the pull back on $C$ of $\pi_\eta$. For all $n$ we have an exact 
sequence:
$$0\to H^0\Cal O_C(D+(n-1)H)\to H^0\Cal O_C(D+nH)\to H^0\Cal O_{Z_\eta}$$
Call $\rho_{\eta,n}$ the last map and put $W_{\eta,n}=$ Im$(\rho_{\eta,n})$, so that 
$\Delta^1_D(n)=\dim W_{\eta,n}$.\par
Choose now $m\gg 0$ and pick a general element $E\in |(m-n)H-D|$. We may identify
$Z_\eta$ and $E$ with their images on $\Gamma$ via $\sigma$.  Surfaces of degree
$m$ through $E$ cut on the general plane $\pi_\eta$ a linear system whose restriction
to $H^0\Cal O_{Z_\eta}$ surjects onto $W_{\eta,n}$. It follows that around a general 
point $\eta$ of
$(\bold P^3)^*$ we may find an \'e{}tale neighbourhood $U$ and a subbundle $\Cal V$
of the $m$-th power of the tautological bundle over the incidence variety in
$\bold P^3\times U$, such that rank$(\Cal V)=\Delta^1_D(n)$ and for any $\eta\in U$ the
fiber $\Cal V_{\eta,n}\subset H^0\Cal O_{\pi_\eta}(n)$ restricts isomorphically to
$W_{\eta,n}$  in $H^0\Cal O_{Z_\eta}$.\par
By assumptions and by construction, $\Cal V_\eta$ satisfies the hypothesis of [CCD2] 
proposition 4, so that, for a general $\eta\in U$,
$$ \Cal V_{\eta,n} \otimes H^0\Cal O_{\pi_\eta}(1) \to H^0\Cal O_{Z_\eta}$$
has rank $\ge \dim W_{\eta,n} +t$.\par
It follows that the restriction map $H^0\Cal O_C(D+(n+1)H)\to H^0\Cal O_{Z_\eta}$
has rank $\ge\dim W_{\eta,n}+t=\Delta^1_D(n)+t$, which gives the claim.
$\diamondsuit$\enddemo

For curves in higher dimensional spaces, we do not have such a precise result
about the descent of the Hilbert function of $D$. Instead, we can use  ideas from
[CCD1] to show the following:

\proclaim{Proposition 2.4}  Let $\Gamma$ be a curve in $\bold P^r$, $\sigma:C\to \Gamma$
its normalization, $D$ a divisor on $C$. Assume
$\Delta^2_D(n)<r-1$ for some fixed $n$; then  either $\Delta^2_D(i)$ vanishes for all 
$i<n$ or it vanishes for all $i>n$.
\endproclaim
\demo{Proof}  As before, consider the exact sequence
$$0\to H^0\Cal O_C(D+(n-1)H)\to H^0\Cal O_C(D+nH)\to H^0\Cal O_{Z}$$
\noindent where $Z$ is the pull back to $C$ of a  general hyperplane section of $\Gamma$; call 
$\rho_n$ the last map 
and put $W_n=$ Im$(\rho_n)$, so that $\Delta^1_D(n)=\dim W_n$. We may identify $Z$ with
its image on $\Gamma$.\par 
Arguing as before, it is not hard to show that, for a general choice of the hyperplane,
the points of $Z$ are in uniform position with respect to $W_n$; it follows by
[CCD1], corollary 1.2, that when $W_{n-1}\neq 0$ and $W_{n}\neq H^0\Cal O_Z$, then
$\dim W_n\ge\dim W_{n-1}+r-1$, i.e. $\Delta^2_D(n)\ge r-1$.\par
On the other hand, observe that $W_{n-1}=0$ implies $h^0\Cal O_C(D+(n-2)H)=h^0\Cal O_C(D+(n-1)H)$ and this is possible only if $\phi_D(n-1)=0$, hence it gives
$\phi_D(i)=0$ for all $i<n$; when $W_n=H^0\Cal O_Z$
then $\Delta^1_D(n)=d$, hence $\Delta^1_D(i)=d$ for all $i\ge d$ and the claim follows.
$\diamondsuit$
\enddemo

\bigskip
\heading 3. Numerical results. \endheading
\smallskip

\proclaim{Definition 3.1} For all $m\in\bold N$ there are uniquely defined numbers
$T(m)$ and $R(m)$ with $0\le R(m)\le T(m)$ such that
$$ m=\frac{T(m)}2 (T(m)+1)+ R(m) = \binom {T(m)+1}2 + R(m).$$
\endproclaim

Let us now turn our attention to non-negative numerical functions  $F:\bold Z\to \bold Z$
which satisfy the following condition \thetag*:

$$\gather F(n)=0\text{ for }|n|\gg 0 \text{ and putting } d=\sum F(i) \text{ we have:}\tag{*}\\
F(n)<t\le s \text{ implies either } \sum_{i=-\infty}^{n-1}F(i)
\le \binom{t-1}2 \text{ or }   \sum_{i=-\infty}^{n-1}F(i)\ge d-\binom t2. 
\endgather$$

The reason for studying these functions is clear from theorem 2.3: $\Delta^2_D$ is
one of them.
We are going to define a function $\Delta$ which minimizes the second sum among the
functions satisfying \thetag*.

\proclaim{Definition 3.2} For given numbers $N$, $s$, $d$ with $d>N\ge 0$ and $d>s^2-s$, define
$T$, $R$, $\rho$, $\epsilon$ as follows:
$$\gather  T=\min(T(N),s-1)\\ R=R(N) \\
     \cases d-1=s(s+1)+\rho s+\epsilon+R-T   &\text{ if }\quad T(N)\le s-1 \\
         d-N-2=\frac s2(s+1)+\rho s+\epsilon &\text{ if } T(N)=s-1\quad \text{ and }\\
                                             &  d-N-1\ge\frac s2(s+1)  \\
          \rho=T(d-N-1)-s\quad\text{ and } &\\
        \epsilon=R(d-N-1)+s-T(d-N-1)-3 &\text{ otherwise }
\endcases \endgather$$  
\noindent where $s>\epsilon\ge 0$.\smallskip\par
Define the function $\Delta(n)=\Delta_{N,s,d}(n)$ as follows:\smallskip

$$  \cases            0 & n<0 \\
                     N+1 & n=0   \endcases  $$
\noindent while for $n$ positive:
$$   \cases
                     T+n+1   &  0<n\le s-T-2 \\
                     s &  s-T-1< n\le s-T+\rho+1 \\
                     2s-(n+T-\rho-1) & s-T+\rho+1< n\le s-T+\rho+\epsilon+1 \\
                     2s-(n+T-\rho) & s-T+\rho+\epsilon+1< n \le 2s-T+\rho \\
                     0 & \text{elsewhere }
\endcases$$
\noindent (observe that some of the previous interval may be empty).\par\smallskip
\endproclaim

When $T(N)\le s-1$, the function $\Delta$ has the following graph:
\vskip3cm

We put here a list of properties of the function $\Delta$ that we are going
to use later. Their proofs are easy consequences of the definition. 

\proclaim{Proposition 3.3}  For all $n>0$: \par
a) $0\le \Delta(n)\le s$.\par
b) $\sum_{i=0}^{+\infty}\Delta(i) = d.$\par
c) $\sum_{i=n}^{+\infty}\Delta(i)>\binom{\Delta(n)}2$ hence 
        $\sum_{i=0}^{n-1}\Delta(i)<d-\binom{\Delta(n)}2$.\par
d) $\sum_{i=0}^{n-1}\Delta(i)\ge\binom{\Delta(n)}2>\binom{\Delta(n)-1}2$.\par
e) For $n\le s-T+\rho$ we have $\sum_{i=1}^n\Delta(i)\le d-\binom s2$.\par
f) If $N'\ge N$ then for all $n\ge 0$: 
     $$\sum_{i=0}^n\Delta_{N',s,d}(i)\ge\sum_{i=0}^n\Delta_{N,s,d}(i).$$\par

\endproclaim

\proclaim{Theorem 3.4} Let $F$ be a non-negative function which satisfies 
condition \thetag*. 
Put $\sum F(i)= d$; fix $N\ge 0$ and assume
$$\sum_{i=-\infty}^{n_0}F(i)\ge \sum_{i=-\infty}^{n_0}\Delta_{N,s,d}(i) \eqno (+)$$
 for some $n_0\ge 0$.  Then for all $n\ge n_0$:
$$\sum_{i=-\infty}^nF(i)\ge\sum_{i=0}^n\Delta(i)$$
\noindent and equality holds if and only if equality holds in \thetag+ 
 and  $F=\Delta$ from $n_0+1$ on.
\endproclaim
\demo{Proof} If the claim does not hold, take $c\ge n_0$ minimal such that
$$\sum_{i=-\infty}^cF(i)<\sum_{i=0}^c\Delta(i)$$
 Necessarily, by the minimality of $c$, we have:
$$ \sum_{i=-\infty}^{c-1}F(i) \ge \sum_{i=0}^{c-1}\Delta(i)\qquad\text{ and }
\qquad F(c)<\Delta(c)\le s,$$
 hence by assumption   we get either 
$$\sum_{i=0}^{c-1}\Delta(i)\le\sum_{i=-\infty}^{c-1}F(i)\le\binom{F(c)}2
\le\binom{\Delta(c)-1}2,$$
\noindent which is excluded by proposition 3.3 d), or 
$$ \sum_{i=-\infty}^{c-1}F(i)\ge d-\binom {F(c)+1}2\ge d-\binom{\Delta(c)}2.$$
Then it follows by 3.3 c):
$$\sum_{i=-\infty}^cF(i)\ge d-\binom {F(c)+1}2+F(c)\ge d-\binom{\Delta(c)}2-
1+\Delta(c)\ge\sum_{i=0}^{c-1}\Delta(i) + \Delta(c)$$
\noindent a contradiction.\par
The second claim is obvious. $\diamondsuit$
\enddemo

\proclaim{Corollary 3.5} Under the previous assumptions, for $n\ge n_0$
$$\sum_{j=-\infty}^n\sum_{i=-\infty}^jF(i)\ge
\sum_{j=0}^n\sum_{i=0}^j\Delta(i).$$
and equality implies $F=\Delta$.
\endproclaim
\demo{Proof} The first assertion is a consequence of theorem 3.4. For the
second, we note that if the equality holds, then certainly $F(n)=\Delta(n)$
for $n> n_0$, by theorem 3.4 and in the interval $(-\infty, 0)$, $\Delta$ has
the smallest second sum among the numerical functions with sum $N+1$.
$\diamondsuit$\enddemo

\bigskip
\heading 4. Applications to linear series. \endheading
\smallskip

In this section, $\Gamma$ is an irreducible, non degenerate curve in $\bold P^3$, $\sigma:C\to
\Gamma$ is its normalization and $h^0\Cal I_\Gamma(s-1)=0$ for a 
fixed number $s$ with $d>s^2-s$.\par
Computing the function $\Delta$ above, one is able to produce several results about 
linear series on $C$. Indeed we have:

\proclaim{Proposition 4.1} Let $\Gamma,C$ be as above, let $D$ be a divisor on $C$, 
assume that $D$ is $H$-primary and $h^0\Cal O_C(D)=N+1$.
Then, for all $n>e(D)$, one has:
$$ \delta\ge\sum_{j=-\infty}^n\sum_{i=-\infty}^j\Delta(i)+g-nd-1,$$
\noindent
with equality if and only if $\Delta^2_D$ coincides with $\Delta_{N,s,d}$.
\endproclaim 
\demo{Proof} The  function $\Delta^2_D$ satisfies condition \thetag* of the previous section,
so if we consider the corresponding function $\Delta:=\Delta_{N,s,d}$, by theorem 3.4
 we know that its second sum bounds from below the Hilbert function of $D$ from 0 on. 
The assertion follows by Riemann-Roch.$\diamondsuit$
\enddemo

In general, the previous proposition leads to (lower) bounds for the degree of $H$-primary  
linear series of given dimension.\par
For instance, if one looks for the gonality of $C$, then one considers linear series $|D|$ 
with $\phi_D(0)=2$, $\phi_D(n)=0$ for $n$ negative.\par\smallskip

In order to apply proposition 4.1, one needs to estimate the number $e(D)$.
It is not hard to find numbers $n$ which satisfy $n>e(D)$: enough to observe, for instance, 
that $e(D)\le e\le (2g-2)/d$.\par 
With more precise information on the level of speciality of $D$, one can expect
better estimates on the degree of $D$. \par
In any event, a reasonable estimate for $e(D)$ is given by the following:

\proclaim{Proposition 4.2} Put $N=h^0\Cal O_C(D)-1$ and define $T,\rho$ as in the previous 
section. Then
$$ h^1\Cal O_C(D+nH)=0 \qquad \text{ for } n\ge 2s-T+\rho-1.$$
\endproclaim
\demo{Proof} For all $n\ge 2s-T+\rho-1$, we have, by theorem 3.4,
$$\Delta^1_D(n)=\sum_{i=-\infty}^n\Delta^2_D(i)\ge
\sum_{i=0}^n\Delta_{N,s,d}(i)=d.$$
The claim follows by proposition 1.4.$\diamondsuit$
\enddemo

\proclaim{Remark 4.3} \rm When $D=0$, we may write the previous function $\Delta_{N,s,d}$, 
starting from $\Delta(0)=1$. It is easy to see that in this case $\Delta$ is the 
function described in [H1], which leads to the computation of the Halphen's bound $G(d;s)$ for the
genus of curves of degree $d>s^2-s$, not lying on any surface of degree $<s$. Recall that

$$G(d;s)=\frac {d^2}{2s}+ \frac {d(s-4)}2+ 1 + \frac {\epsilon}2 (s-\epsilon-1+\frac \epsilon s)$$

\noindent where $d=ks-\epsilon$, with $0\leq \epsilon\leq s-1$. \par In
other words, by applying propositions 4.1 and 4.2 above, for $D=0$, one finds the usual
Halphen's bound for the genus of space curves not contained  on surfaces of degree $<s$.
\endproclaim

Using proposition 4.2, and computing the second sum of $\Delta$, up to $2s-T+\rho$, one gets:

\proclaim{Theorem 4.4} Any $H$-primary linear series $g^r_\delta$ on $C$ satisfies:
$$\multline \delta \ge
g+(T-\rho-2s)d-1+(\rho+1)\big[R-T+\frac{3s^2+s}2+\frac s2(\rho+2)\big]- \\
-\frac{(s-\epsilon)(s-\epsilon-1)}2+ s^3+\frac32s^2+\frac s2-2sT+2sR-RT-\frac{T^3-3T^2+2T}6
\endmultline$$
\noindent where $R,T$ are defined as above.
\endproclaim

This bound is sharp, as we shall see soon. 
\par\smallskip

Since any linear series of dimension 1 or 2 is $H-$primary, one gets:

\proclaim{Corollary 4.5} Let $\Gamma$ be an irreducible space curve of degree $d$ and
geometric genus $g$; assume that it does not lie on surfaces of degree $<s$ and $d>s^2-s$. 
Call $C$ a normalization of $\Gamma$.\par
If $D$ is a  divisor of degree $\delta$ on $C$, defining a  linear
series $g^1_{\delta}$, then $N=1$, $T=1$, $R=0$ and:
$$\delta\ge g-s^3+\frac 32s^2+\frac32s-\frac32\rho s^2-s\epsilon +2\rho s+\frac\epsilon2
-\frac{\epsilon^2}2-\frac s2\rho^2-\rho\epsilon-\rho-2.$$
If $D$ is a  divisor of degree $\delta$ on $C$, defining a linear
series $g^2_{\delta}$, then $N=2$, $T=1$, $R=1$ and:
$$\delta\ge g-s^3+\frac 32s^2+\frac32s-\frac32\rho s^2-s\epsilon +2\rho s+\frac\epsilon2
-\frac{\epsilon^2}2-\frac s2\rho^2-\rho\epsilon-\rho-1$$
\noindent (clearly $\epsilon$ and $\rho$ differ in the two expressions).\endproclaim 

Let us discuss some examples.

\proclaim{Example 4.6} \rm When $C=\Gamma$ is a smooth complete intersection of surfaces of 
degree $s,p$, with $p\ge s (\ge 2)$,  we find from the previous expression that, if $|D|$ is 
a $g^1_\delta$ on $C$, then
$$\delta \ge sp-p. $$
\noindent This is exactly the bound found in [B] (see also [CL]  and [L], pag.23). \par
This bound is attained when $C$ has a $p$-secant line $\ell$ and $D$ is cut out by the
planes through $\ell$, off the points cut out  by $\ell$ on $C$. In fact,
this is the unique kind of $g^1_{sp-p}$ on $C$. Indeed, since $D$ satisfies
$\Delta^2_D(n)=\Delta(n)$ for all $n$,  by 1.5 one computes
$h^1\Cal O_C(D+(s+p-5)H)>0$. So $\omega_C\otimes\Cal O_C(-(s+p-5)H-D)\simeq\Cal O_C(H-D)$ is 
effective and clearly $H-D$ is contained in 2 independent planes.\par
Similarly, for a $g^2_\delta$ one has:
$$\delta\ge sp-1$$
\noindent and the minimum is attained when the divisors lie on a plane, i.e. when the
series is cut by the planes through a point.
\endproclaim

Next, we apply 4.4 to {\it Halphen curves.}
 
\proclaim{Theorem 4.7}  Let  $C=\Gamma$ be a smooth {\rm Halphen curve} of genus $g=G(d;s)$,
that is, a curve of degree $d>s^2-s$ and maximal genus among those not lying on surfaces 
of degree $<s$. Put $d-1=ms+e$, $0\le e<s$ and assume that $C$ is not complete intersection, 
i.e. $e\neq s-1$.\par
Then any linear series $g^1_\delta$ on $C$ satisfies $\delta\ge ms-m$ and
this bound is sharp: it is attained if and only if $C$ has a $(m+e+1)$-secant line $\ell$, 
by the series  $|H-C\cdot \ell|$, where $H$ is cut out on $C$ by a plane.\par
For any linear series $g^2_\delta$, we have $\delta\ge d-1$. Also this bound is sharp;
if either $s>3$ and $m>2$ or $e<s-2$, then every $g^2_{d-1}$ on $C$ is cut out by the planes
through a fixed point. \endproclaim
\demo{Proof} We know that $C$ is directly linked to a plane
curve $C'$ of degree $s-e-1$ by surfaces of degree $s$ and $m+1$ (se [H1]).\par 
The two bounds are a straightforward consequence of theorem 4.4; for the sharpness,
one has only to observe that there are Halphen curves with a $(m+e+1)$-secant line: they
can be constructed as the residual of a plane curve of degree $s-e+1$ containing a line, 
in a complete intersection of type $(s,m+1)$.\par
It remains only to classify series of minimal degree.\par
Assume that a divisor $D$ defines a $g^1_{ms-m}$ on  a smooth Halphen curve $C$; 
then we only need to prove that $H-D$ is effective, for then it must be contained on a line, 
because $C$ is arithmetically normal. \par
The canonical series is cut on $C$ by surfaces of 
degree $s+m-3$ passing through $C'$ and, by theorem 3.4, we know the function
$\Delta_D^2$, which must coincide with the corresponding $\Delta_{1,s,d}$. 
Put $E=C\cdot C'$ and consider the following cases:\par
\medskip
\noindent  (i) case $e=s-2$; as in the previous example, computing $\Delta^2_D$ one finds that
$h^1\Cal O_C(D+(s+m-4)H)>0$, which, by Serre duality, means that $H-D-E$ is effective;\par 
\medskip
\noindent (ii) case $e=s-3$; arguing as above, one gets here $h^0\Cal O_C(2H-D-E)\ge 2$ and since
$C'$ has degree 2, it follows that there is a quadric containing $D$ and containing
the plane of $C'$; so $D$ must lie on a plane;\par
\medskip
\noindent (iii) case $e<s-3$; in this case $h^1\Cal O_C(D+(s+m-5)H)>0$, so $2H-D-E$ is effective; 
since deg($C')>2$, every quadric through $E$ contains the plane of $C'$, hence $H-D$ 
is effective.\par
\medskip
This complete the proof of the claim for the $g^1_\delta$'s.\par
Take now a divisor $D$ which defines a $g^2_{d-1}$ on $C$. As above, since $\Delta^2_D=
\Delta_{2,s,d}$, one finds that $h^0\Cal O_C(2H-D-E)>0$, hence there is an effective
divisor in $|2H-D-E|$ which we denote again by $2H-D-E$. Consider the following cases:\par
\medskip 
\noindent (i) case $e<s-3$; here $C'$ has degee $>2$ hence any quadric through
$D+E$ contains the plane of $C'$, so $H-D$ is effective and we are done;\par
\medskip
\noindent (ii) case $e=s-3$;  $C'$ is a conic and the effective divisor $2H-D-E$ imposes 2
conditions on $|2H-E|\supset (H-E)+|H|$. If $2H-D-E$ has no common points with $H-E$, then it
imposes 2 conditions on $|H|$, hence it must span a line, which is impossible for degree
reasons. If $2H-D-E$ has a common point with $H-E$, then any  
quadric through $D+E$  contains the plane of $C'$, hence
$H-D$ is effective;\par
\medskip
\noindent (iii) case $e=s-2$; then $C'$ is a line and $h^0\Cal O_C(2H-E)=7$ and one computes 
$h^0\Cal O_C(2H-D-E)=2$. Assume that $|2H-D-E|$ is base  point free; then it cannot impose 4
conditions to $H$, otherwise by Castelnuovo's lemma (see [C1]) it would impose at least 5
conditions to $2H-E$, a contradiction, since this would imply $h^0\Cal O_C(D)\le 2$; on the other
hand, if  $2H-D-E$ lies in some plane, then $|2H-D-E|$ is cut by the planes through a fixed 
 $(m+s-2)$-secant line; when $s>3$ and $m>2$ this line cannot exists, because
$C'$ is $(m+s-1)$-secant and $C\cup C'$ is complete intersection of type $m+1, s$.
Hence $|2H-E-D|$ has a fixed point $P$, i.e. $|2H-E-D|=P+|F|$, where $|F|$ is a $g^1_{ms-m}$;  
as we saw above,  $|F|$ must be cut by the planes through a  $(m+s-1)$-secant line, which  
must be $C'$ when $s>2$; it follows that $H-D$ is effective.$\diamondsuit$
\enddemo

Notice that the hypothesis $s>3$ and $m>2$, i.e. the hypothesis that the curve does not lie on
a cubic surface, in the statement of theorem 4.7
 above is necessary: for instance, a smooth curve of genus 2 and degree 5 is a Halphen curve,
the complete intersection of a quadric and a cubic surface off a line, and it has a
2-dimensional family of $g^2_4$ on it, so the general such series is not cut out by planes
through a point of the curve. In general, any smooth Halphen
curve $C$ of degree $d=3m+2$ on a cubic surface has linear series $g^2_{d-1}$ which are not cut
out by planes: namely, as the proof of theorem 4.7 shows, any series of the type
$|H-E+E'|$, where $E$ is cut out on $C$ by the $(m+2)-$secant line and $E'$ is cut out on $C$
by another line on the cubic.  \par

Halphen curves  are not the only
curves which attain the bound of corollary 4.4. The curves  in next example  are not even
arithmetically normal.

\proclaim{Example 4.8}\rm Let $C=\Gamma$ be a smooth curve which is directly linked to a pair 
$Y$ of disjoint lines, by surfaces of degree $s,p$, with $2< s\le p$. Then we have $d=sp-2$, 
$g=-1+(s+p-4)(sp-4)/2 $ and $C$ is not contained on surfaces of degree $<s$.
For a $g^1_\delta$ on $C$ one must have, by corollary 4.2, $\delta\ge ps-p-s$.\par
On the other hand, it is not hard    to compute that both  lines of $Y$
are $(s+p-2)$-secant to $C$, hence the set of planes through any of them
cuts on $C$ a $g^1_{ps-p-s}$. \par
Arguing as in the proof of theorem 4.7, it turns out that these series are the unique pencils 
of minimal degree on $C$. We omit the details.
\endproclaim

For curves of relatively low genus with respect to the degree, the estimate of $e(D)$ given 
by proposition 4.2 is not sharp. This affects the sharpness of the bounds
given in 4.4 and 4.5 too; so in order to improve them, one has to improve
the estimate on $e(D)$.\par
Let us show for example what happens for curves on a smooth quadric.

\proclaim{Example 4.9}  \rm Let $Q$ be a smooth quadric surface and let $C$ be a
curve of type $(a,b)$ on $Q$, with $a\le b$. Then we know that $d=a+b$,
$g=ab-a-b+1$, $e(0)=a-2$ and $s=2$. If $D$ defines a $g^1_\delta$ on $C$, the corresponding
function $\Delta$ can be easily computed: it is a sequence of 2's from 0 to $\big[\frac d2\big]$
with a final 1 if $d$ is odd, and 0 elsewhere.\par
Using corollary 4.5 one computes
$\delta\ge a$, which is the expected result.\par
Since we know $\omega_C$, we can prove as above that for $a<b$,
a series of minimal degree must be cut by a ruling of $Q$.\par
For a $g^2_\delta$, as above, one shows similarly  that, for $a<b$, then $\delta\ge 2a$, 
with the minimum attained only by the series cut by pair of lines in a ruling; when $a=b$, the 
bound 
is $\delta\ge d-1$ and, as above, the minimum is attained by the series cut by the planes 
through a point.
\endproclaim

\proclaim{Remark 4.10} \rm The bounds that one finds for $\delta$ in the case
$e(D)\ge s-T-1$ are strictly increasing with $s$. It follows that, for fixed
genus $g$ and degree $d>s^2-s$ of $\Gamma$, the linear
series attaining these bounds must live on curves lying on surfaces of degree $s$.
\endproclaim

\proclaim{Example 4.11}  \rm Using our bound for
computing the minimal degree of a $H$-primary $g^3_\delta$ on a Halphen curve $C$, 
one gets:
$$\delta   \ge 2d-6,  $$ which is still sharp: look at the series cut by quadrics through 6 
general points of $C$. In particular, when $d>6$,  every $g^3_d$ on $C$ cannot be $H$-primary.
\endproclaim

This last fact leads to the following result, which should be compared with theorem 3.1 from
[CL]:

\proclaim{Proposition 4.12} If $C$ is a smooth Halphen curve of degree $d>6$, i.e.
a curve of maximal genus among those not lying on surfaces of degree $<s$, for some 
$s$ with $d>s^2-s$, then $|H|$ is the unique $g^3_d$ on $C$.\par
In other words, there is at most one way to embed a curve in $\bold P^3$ as a
Halphen curve of degree $d>6$.\endproclaim

When the linear series $g^N_\delta$ is not $H$-primary, one gets lower bounds for
$\delta$ using corollary 3.5.

\proclaim{Proposition 4.13}  Let $D$ define a $g^N_\delta$ on $C$ and let 
$$Q = \max\{n: N+1\ge \sum_{j=0}^n\sum_{i=0}^j\Delta_{0,s,d}(i)\}.$$
Then for all $n>e(D)$:
$$\delta \ge g-nd -1 +\sum_{j=0}^{Q+n}\sum_{i=0}^j\Delta_{0,s,d}(i) $$
\noindent with equality if and only if $N+1=\sum_{j=0}^Q\sum_{i=0}^j\Delta_{0,s,d}(i)$ 
and $\Delta^2_D(i)=\Delta_{0,s,d}(i-Q)$.\par In particular we get:
$$\delta\ge g-(2s+\rho-Q-1)d -1 +\sum_{j=0}^{2s+\rho-1}\sum_{i=0}^j\Delta_{0,s,d}(i) $$
\endproclaim
\demo{Proof} By our assumptions, there exists some number $c$, with 
$-Q\le c\le 0$, such that:
$$\Delta^1_D(c)\ge \sum_{i=0}^{c+Q}\Delta_{0,s,d}(i).$$
If we define  $F(i)=\Delta^2_D(i-Q)$, then 
$$\sum^{c+Q}_{i=-\infty}F(i)\ge \sum_{i=0}^{c+Q}\Delta_{0,s,d}(i)$$
and the first assertions follows by  corollary 3.5.\par
For the second, just note that, by construction,
$$\sum_{i=0}^{2s+\rho-1}\Delta_{0,s,d}(i)=d. \diamondsuit$$
\enddemo

Computing this bound for the case of a $g^3_\delta$ on a Halphen curve, one gets $\delta\ge d$,
so the bound is sharp, but we do not expect sharpness in more general situations.

\proclaim{Remark 4.14} \rm It is possible to analyse other classes of space curves
with the method we used here, but the computations involved become often
quite cumbersome.\par
For example, one can consider smooth arithmetically Cohen Macaulay curves with given arithmetical
characters of the Hilbert-Burch resolution of their coordinate ring. In this case one can use the
bound on $e$ given in [CGO], in order to get bounds for the degree of primary linear series. Once
these computations have been performed, the problem of sharpness remains. However, at least for
linear series of dimension $N\le 3$, we expect the answer to this problem to be {\it the natural
one}. For example, for $N=3$, we already know, from [CL], that in general there is  no non primary
$g^3_d$. Similarly, we expect that curves with minimal gonality should in general be obtained by a
direct linkage to a curve containing a line.
\endproclaim

\medskip

\heading  5. Results in $\bold P^r$. \endheading

When the curve sits in $\bold P^r$  we can obtain, using proposition 2.4,
a  bound for the degree of a linear series of given dimension. 
Since we are unable to prove a  precise descent lemma like
theorem 2.3, the bound we obtain will  be sharp essentially only for
curves lying on minimal scrolls.\par\smallskip

The non-negative numerical functions $F$ that we consider now are those satisfying the following
condition:
$$ \text{If }  F(n)<r-1 \text{ then either }F(m)=0\quad\forall m<n \text{ or }F(m)=0\quad\forall 
    m>n\tag{**}$$

We define here the function $\Delta=\Delta_{N,d}^r$, for $d\ge r+N$, by setting:

$$\gather d-N-2=\rho(r-1)+\epsilon\qquad 0\le\epsilon<r-1 \\
\Delta(n)=\cases N+1  & n=0 \\
                   r-1 &  1\le n\le \rho \\
                   \epsilon+1 & n=\rho+1 \\
                   0 & elsewhere
\endcases \endgather$$
We are able to reproduce easily all the previous argument, obtaining:

\proclaim{Theorem 5.1} Let $F$ be a non-negative function which satisfies condition \thetag{**}.
Put $\sum F(i)= d$; fix $N>0$ and assume 
$$\sum_{i=-\infty}^{n_0}F(i)\ge\sum_{i=0}^{n_0}\Delta_{N,d}^r(i) \eqno (++)$$
for some $n_0\ge 0$. Put $\Delta=\Delta^r_{N,d}$.
Then for all $n\ge n_0$:
$$\sum_{i=-\infty}^nF(i)\ge\sum_{i=0}^n\Delta(i)\tag i$$
\noindent and equality holds if and only if equality holds in \thetag {++}
and $F=\Delta$ from $n_0$ on.
$$\sum_{j=-\infty}^n\sum_{i=-\infty}^jF(i)\ge\sum_{j=0}^n\sum_{i=0}^j\Delta(i)\tag {ii}$$
and equality implies $F=\Delta$.\par
Moreover, assume that, for some $n_0\ge 0$ one has
$$\sum_{j=-\infty}^{n_0}\sum_{i=-\infty}^jF(i)=\sum_{j=0}^{n_0}\sum_{i=0}^j\Delta_{0,s,d}(i).$$
Then for all $n\ge n_0$:
$$\sum_{j=-\infty}^n\sum_{i=-\infty}^jF(i)\ge\sum_{j=0}^n\sum_{i=0}^j\Delta_{0,s,d}(i)
\tag{iii}$$
with equality if and only if $F=\Delta$.\par
Finally, we have $e(D)<\rho$.
\endproclaim

Turning to our curves $C$, we know by proposition 2.4 that $\Delta^2_C$ satisfies condition 
\thetag{**}. By Riemann-Roch we have:

\proclaim{Corollary 5.2}  Assume that $D$ defines a $H$-primary $g^N_\delta$ on the normalization $C$
of a non-degenerate curve $\Gamma\subset\bold P^r$. Then for all $n>e(D)$:
$$ \delta\ge\sum_{j=-\infty}^n\sum_{i=-\infty}^j\Delta(i)+g-nd-1$$
\noindent with equality if and only if $\Delta^2_D=\Delta$.
In particular, we have:
$$ \delta\ge g-\epsilon\rho+N-\rho-(\rho^2-\rho)\frac{r-1}2 $$
\endproclaim

\proclaim{Example 5.3} \rm If we apply corollary 5.2 to find the minimal degree
for $H$-primary linear series $g^N_\delta$ on Castelnuovo curves in $\bold P^r$, then we find in 
many cases the expected value (see [A] for the case of pencils).\par
Indeed we recall (see [H2]) that such a curve is smooth and lies on a surface $S$ of minimal 
degree $r-1$; when
$r>5$ and $S$ is smooth, then it is a scroll whose Picard group is generated by two classes $h,f$,
where $f$ represents the ruling and $h$ the hyperplane class. In this case,  
Castelnuovo curves on $S$ are in the classes $ah-(r-2+b)f$, with $b = 0,\dots,r-2  $ and for 
a $H$-primary $g^N_\delta$, the minimal value of
$\delta$ is computed, from the previous formula; one gets:
$$ \delta\ge \cases aN & b\ge N \\ (a-1)N+b & b<N \endcases$$
\noindent When $N=1$, this is the value found in [A]; it is attained by the linear series cut 
by $f$. Furthermore, since we know the canonical class of $S$,  $\omega_S=-2h+(r-3)f$, we 
may also compute that a linear pencil of minimal degree, whose function $\Delta_D^2$ is 
then equal to $\Delta$, must be cut by the ruling. This is done as for the examples considered in
\S 4.\par
For linear series of dimension $N=2,\dots,r-1$ (necessarily
they are $H$-primary)  we find that the minimal degree 
 is achieved when the series is a multiple of the ruling, in the case $b\ge N$,
while for $b<N$ the bound is always attained, for $b=0$, on curves lying on a cone, by 
the $N$-th multiple of the linear series cut out by the ruling, off the vertex of the cone.\par
Furthermore, one can prove directly that the previous bound is always attained, for Castelnuovo 
curves and for all $N=2,...,r-1$, when $r\le 5$ (recall that in  $\bold P^5$ one may also
consider curves lying on a Veronese surface).
 On the other hand a smooth Castelnuovo curve $C$ in $\bold P^6$, of degree 12 has genus 7 and 
the previous formula  tells us that  a $g^2_{\delta}$ on $C$ has degree $\delta\ge 5$. This is 
certainly true, but the bound cannot be sharp for there are no  $g^2_5$ on a smooth curve of
 genus 7, unless the curve is hyperelliptic, which is not our case.\par
The reason for this relies in the fact that our method uses a reduction to
some hyperplane section of $C$, hence it is not sensitive to the isolated singularities
of the curve; in the previous example, if we allow a node for $C$, then the genus becomes 6 and 
the formula tells $\delta\ge 4$ for a $g^2_\delta$; this bound is
attained on the hyperelliptic curve of degree 12 which is residual
of 3 rulings on a quintic cone. Notice that this curve has maximal arithmetic
genus among non-degenerate curves in $\bold P^6$.
\endproclaim

\proclaim{Example 5.4} \rm When the curve $C$ lies on a minimal scroll but it is not
a curve of maximal genus, then the previous formulae are not sharp. One has to
put into the picture the number $e(D)$, which is easily computed:
 we get $e(D)=a-2$ for a curve of type $ah+bf$, $b\ge 0$.\par
Using the function $\Delta$, one gets that the minimal gonality, for such a curve, is exactly $a$.
\endproclaim

\proclaim{Example 5.5} \rm When $N\ge r$, we have to distinguish the $H$-primary and the non 
$H$-primary cases and use theorem 5.1 \thetag{iii} for this last one. 
For instance, when $C$ is a Castelnuovo curve, as in example 5.3, for $N=r$ one finds
that the minimal degree for a $g^r_\delta$  on $C$ is $\delta=d$ and the unique series of 
smallest degree is the hyperplane series. Indeed for $H$-primary series we get
$\delta \ge ra-r+b $, which is bigger than $d$. 
\endproclaim

\end